\documentclass[12pt,a4paper,reqno]{amsart}
\usepackage[utf8]{inputenc}
\usepackage{enumitem}
\usepackage{amsmath}
\usepackage{amsfonts}
\usepackage{amssymb}
\usepackage[all]{xy}
\usepackage{xfrac}
\usepackage[yyyymmdd,hhmmss]{datetime}
\usepackage{centernot}
\usepackage[english]{babel}
\usepackage{dsfont}
\usepackage{bbold}
\usepackage{soul}
\usepackage{lmodern}

\newtheorem{theorem}{Theorem}[section]

\newtheorem{lemma}[theorem]{Lemma}

\newtheorem{proposition}[theorem]{Proposition}

\newcommand{\e}{\wedge}
\newcommand{\ou}{\vee}

\usepackage{mathtools}
\newcommand{\eq}[1]{\stackrel{\mathclap{\normalfont\mbox{{\tiny #1}}}}{=}}
\newcommand{\impli}[1]{\stackrel{\mathclap{\normalfont\mbox{{\tiny #1}}}}{\Rightarrow}}

\author{J. P. Fatelo and N. Martins-Ferreira}

\address{School of Technology and Management, Centre for Rapid and Sustainable Product Development - CDRSP, Polytechnic Institute of Leiria, P-2411-901 Leiria, Portugal.}

\email{martins.ferreira@ipleiria.pt}

\title[]{A new look at Ternary Boolean Algebras}

\subjclass[2020]{Primary 06E05, 
                          06D30; 
								  Secondary 03G25,  
								            03C48}  

\keywords{Boolean algebra, non-symmetric ternary operation, ternary Boolean algebra, algebraic structure.}

\thanks{  This work is supported by Fundação para a Ciência e a Tecnologia (FCT\-UID-Multi-04044-2019), Centro2020 (PAMI -- ROTEIRO\-/0328\-/2013\-- 022158) and Polytechnic of Leiria through the projects CENTRO\--01\--0247-FEDER: 069665, 069603, 039958, 039969, 039863, 024533 and also ESTG and CDRSP}
 
\begin{document}

\begin{abstract}
We present a new approach to ternary Boolean algebras in which negation is derived from the ternary operation. The key aspect is the replacement of 
complete commutativity by other axioms that do not require the ternary operation to be symmetric.

\end{abstract}

\maketitle

\section{Introduction}

In 1947, Grau \cite{Grau} considered an algebraic structure $(A,p,\bar{()})$ consisting of a set $A$, a ternary operation $p\colon{A^3\to A}$ and a unary operation $\bar{()}\colon{A\to A}$ satisfying the following axioms
\begin{enumerate}[label={\bf (A\arabic*)}]
\item\label{A1} $p(a,b,p(c,d,e))=p(p(a,b,c),d,p(a,b,e))$
\item\label{A2}
$p(a,b,b)=p(b,b,a)=b$
\item\label{A3}
$p(a,b,\bar{b})=p(\bar{b},b,a)$.
\end{enumerate}
This structure was called ternary Boolean algebra. For every chosen element $0\in A$ with $a\wedge b=p(a,0,b)$ and $a\vee b=p(a,\bar{0},b)$, the system $(A,\vee,\wedge,\bar{()},0,\bar{0})$ is a Boolean algebra. Every other choice of an element $0\in A$ gives an isomorphic Boolean algebra. Furthermore, the following condition (called complete commutativity) is derived from the axioms
\begin{equation}\label{eq:complete commutativity}
p(a,b,c)=p(a,c,b)=p(c,a,b)
\end{equation}
and the ternary operation is uniquely determined as
\begin{eqnarray}\label{eq:p(a,b,c) as Grau formula}
p(a,b,c)=(a\wedge b)\vee (b\wedge c) \vee (c\wedge a).
\end{eqnarray}

In the same year, Birkhoff and Kiss \cite{Birk} considered a system $(A,p,0,1)$ satisfying the axioms
\begin{enumerate}[label={\bf (B\arabic*)}]
\item\label{B1} $p(0,a,1)=a$
\item\label{B2} $p(a,b,a)=a$
\item\label{B3} $p(p(a,b,c),d,e)=p(p(a,d,e),b,p(c,d,e))$
\item\label{B4} $p(a,b,c)=p(b,a,c)=p(b,c,a)$
\end{enumerate}
and showed that it is the same as a distributive lattice. {The meet and join operations are defined again as $a\wedge b=p(a,0,b)$ and $a\vee b=p(a,1,b)$ and the distributive lattice is bounded with $0$ and $1$ as bottom and top elements, respectively. It is a Boolean algebra when there exists a unary operation $\bar{()}\colon{A\to A}$ such that $p(a,0,\bar{a})=0$ and $p(a,1,\bar{a})=1$. Once more, the ternary operation $p$ is determined by the formula (\ref{eq:p(a,b,c) as Grau formula}).
Note that \ref{A1} and \ref{B3} are equivalent under \ref{B4}. Birkhoff and Kiss called symmetry to axiom \ref{B4} which is the same as Grau's complete commutativity (\ref{eq:complete commutativity}).

The purpose of this note is to show that condition~\ref{B4} can be replaced by another type of conditions, namely $p(a,0,b)=a$ and $p(a,1,b)=b$.  In this new setting, the meet and join operations are defined differently as:
\begin{equation*}\label{eq:meet-join}
a\wedge b=p(0,a,b)\quad \textrm{and}\quad a\vee b=p(a,b,1). 
\end{equation*}
Since complete commutativity is no longer present, conditions \ref{A1} and \ref{B3} must be carefully reformulated as
\begin{equation*}\label{eq:homo}
p(a,p(b_1,b_2,b_3),c)=p(p(a,b_1,c),b_2,p(a,b_3,c)).
\end{equation*}
This gives rise to a different formula for the ternary operation $p$ with the advantage that negation may be derived as
\begin{equation*}\label{eq:bar}
\bar{a}=p(1,a,0).
\end{equation*}
The resulting formula for the ternary operation is then
\begin{equation*}\label{eq:our formula for p}
p(a,b,c)=(\bar{b}\wedge a)\vee (b\wedge c).
\end{equation*}

In summary, our result is the following theorem.
\begin{theorem}\label{thm:1}
Let $(A,p,0,1)$ be a system consisting of a set $A$, together with a ternary operation $p$ and two constants $0,1\in A$ satisfying the following conditions:
\begin{enumerate}[label={\bf (C\arabic*)}]
\item\label{C1}\label{smA2} $p(0,a,1)=a$
\item\label{C2}\label{smA3} $p(a,b,a)=a$
\item\label{C3}\label{smA7} $p(a,p(b_1,b_2,b_3),c)=
p(p(a,b_1,c),b_2,p(a,b_3,c))$
\item\label{C4}\label{smA4}\label{smA5} $p(a,0,b)=a=p(b,1,a)$.
\end{enumerate}
For $a+b=p(a,b,\bar a)$, $a\wedge b=p(0,a,b)$, $a\vee b=p(a,b,1)$, \mbox{$\bar{a}=p(1,a,0)$},
the following conditions are equivalent:
\begin{enumerate}[label={\rm(\roman*)}]
\item\label{T5} The system $(A,+,\wedge,0,1)$ is a Boolean ring
\item\label{T1} The system $(A,\vee,\wedge,\bar{()},0,1)$ is a Boolean algebra
\item\label{T2} $p(a,b,c)=(\bar{b}\wedge a)\vee (b\wedge c)$
\item\label{T3} $p(0,a,b)=p(a,a,b)$
\item\label{T4} $p(a,b,1)=p(a,b,b)$.
\end{enumerate}
\end{theorem}
\noindent 
The proof of this theorem is detailed in the next section.

\section{Proof of Theorem \ref{thm:1}}\label{section:proof}
We start with some properties derived from conditions \ref{C1} to~\ref{C4} of Theorem \ref{thm:1}. The following notation is used:
$$\overline{a}=p(1,a,0),\quad a\wedge b= p(0,a,b),\quad a\vee b=p(a,b,1),\quad a+ b=p(a,b,\bar a).$$
\begin{proposition}\label{lemma:1}
If $(A,p,0,1)$ verifies conditions \ref{C1} to~\ref{C4} then:
\begin{eqnarray}
\label{L1}\bar{1}=0 &,& \bar{0}=1\\
\label{L2} \overline{\overline{a}}&=&a\\
\label{L3} p(c,b,a)&=&p(a,\bar{b},c)\\
\label{L4} \overline{p(a,b,c)}&=&p(\overline{a},b,\overline{c})\\
\label{L5} \overline{p(a,b,c)}&=&p(\overline{c},\overline{b},\overline{a})\\
\label{L6} \overline{a\wedge b}=\overline{b}\vee \overline{a}&,&\overline{a\vee b}=\overline{b}\wedge \overline{a}\\
\label{L7} (A,\wedge,1)&\textrm{is}&\textrm{a monoid}\\
\label{L8} (A,\vee,0)&\textrm{is}&\textrm{a monoid}\\
\label{L9} a\wedge 0=&0&=0\wedge a\\
\label{L10} a\vee 1=&1&=1\vee a\\
\label{monoid-plus}(A,+,0)&\textrm{is}& \textrm{a monoid}.
\end{eqnarray}
\end{proposition}
\begin{proof}
In each step of the proof, the needed property when required is written above the corresponding equality.
$$\begin{array}{l}
\bar{1}=p(1,1,0)\ \eq{\ref{C4}}\ 0,\quad \bar{0}=p(1,0,0)\ \eq{\ref{C4}}\ 1\\[5pt]
\overline{\overline{a}}=p(1,p(1,a,0),0)\ \eq{\ref{C3}}\ p(p(1,1,0),a,p(1,0,0))\ \eq{\ref{C4}}\ p(0,a,1)\ \eq{\ref{C1}}\ a\\[5pt]
p(a,\overline{b},c)=p(a,p(1,b,0),c)\ \eq{\ref{C3}}\ p(p(a,1,c),b,p(a,0,c))\ \eq{\ref{C4}}\ p(c,b,a)\\[5pt]
\overline{p(a,b,c)}=p(1,p(a,b,c),0)\ \eq{\ref{C3}}\ p(p(1,a,0),b,p(1,c,0))=p(\overline{a},b,\overline{c}).
\end{array}$$
Property (\ref{L5}) is just a combination of (\ref{L3}) and (\ref{L4}). Next is the proof of Properties (\ref{L6}):
$$\begin{array}{l}
\overline{a\wedge b}=\overline{p(0,a,b)}\ \eq{(\ref{L5})}\ p(\bar{b},\bar{a},\bar{0})\ \eq{(\ref{L1})}\ p(\bar{b},\bar{a},1)
=\bar{b}\vee\bar{a} \\[5pt]
\overline{a\vee b}=\overline{p(a,b,1)}\ \eq{(\ref{L5})}\ p(\bar{1},\bar{b},\bar{a})\ \eq{(\ref{L1})}\ p(0,\bar{b},\bar{a})
=\bar{b}\wedge\bar{a}.
\end{array}$$
With respect to (\ref{L7}) and (\ref{L8}), we have associativity
$$\begin{array}{rcccl}
(a\wedge b)\wedge c&=&p(0,p(0,a,b),c)&\eq{\ref{C3}}& p(p(0,0,c),a,p(0,b,c))\\[5pt]
                   &\eq{\ref{C4}}&p(0,a,p(0,b,c))&=&a\wedge(b\wedge c)\\[5pt]
a\vee(b\vee c)&=&p(a,p(b,c,1),1)&\eq{\ref{C3}}& p(p(a,b,1),c,p(a,1,1))\\[5pt]
              &\eq{\ref{C4}}&p(p(a,b,1),c,1)&=&(a\vee b)\vee c,
\end{array}$$
and identities
$$\begin{array}{l}
a\wedge 1=p(0,a,1)\ \eq{\ref{C1}}\ a,\quad 1\wedge a=p(0,1,a)\ \eq{\ref{C4}}\ a\\[5pt]
a\vee 0=p(a,0,1)\ \eq{\ref{C4}}\ a,\quad 0\vee a=p(0,a,1)\ \eq{\ref{C1}}\ a.
\end{array}$$
For properties (\ref{L9}) and (\ref{L10}), the proof is:
$$\begin{array}{l}
a\wedge 0=p(0,a,0)\ \eq{\ref{C2}}\ 0,\quad 0\wedge a=p(0,0,a)\ \eq{\ref{C4}}\ 0\\[5pt]
a\vee 1=p(a,1,1)\ \eq{\ref{C4}}\ 1,\quad 1\vee a=p(1,a,1)\ \eq{\ref{C2}}\ 1.
\end{array}$$
The structure $(A,+,0)$ is a monoid since:
$$a+0=p(a,0,\bar a)\ \eq{\ref{C4}}\ a,\quad
0+a=p(0,a,1)\ \eq{\ref{C1}}\ a$$
\begin{eqnarray*}
(a+b)+c&=&p(p(a,b,\bar a),c,\overline{p(a,b,\bar a})\ \eq{(\ref{L5})}\ p(p(a,b,\bar a),c,p(a,\bar b,\bar a))\\
			 &\eq{\ref{C3}}& p(a,p(b,c,\bar b),\bar a)=a+(b+c).
\end{eqnarray*}
\end{proof}


In addition to the structure defined by conditions \ref{C1} to~\ref{C4}, in the following lemma \mbox{$\bar a=p(1,a,0)$} is assumed to be the Boolean complement of $a\in A$. When this is the case $\e$ and~$\ou$ are idempotent.
Note that the de Morgan's laws (\ref{L6}) imply a duality between $\wedge$ and $\vee$. In particular, 
$\wedge$ is idempotent if and only if $\vee$ is idempotent too. 
\begin{lemma}\label{lemma:complement}
Let $(A,p,0,1)$ be a system verifying conditions \ref{C1} to~\ref{C4}. If every $a\in A$ verifies the relations
\begin{equation}\label{complement}
\bar a\e a=p(0,\bar a,a)=0\quad \textrm{and}\quad \bar a\ou a=p(\bar a,a,1)=1,
\end{equation}
then idempotency holds:
\begin{equation}\label{idem}
p(0,a,a)=a\wedge a=a\quad \textrm{and}\quad p(a,a,1)=a\vee a=a.
\end{equation}
\end{lemma}
\begin{proof}
$$\begin{array}{rcccl}
a\wedge a&=&p(0,a,a)\quad \eq{(\ref{complement}),\ref{C4}}\quad p(p(0,\bar a,a),a,p(0,1,a))\\[5pt]
         &\eq{\ref{C3}}&\hspace{-23pt}p(0,p(\bar a,a,1),a)\ \eq{(\ref{complement})}\ p(0,1,a)\ \eq{\ref{C4}}\ a
\end{array}$$
Idempotency of $\vee$ is obtained similarly or using (\ref{L6}).
\end{proof}

The next lemma shows how the condition $\ref{T3}$ of Theorem~\ref{thm:1} turns the structure of axioms \ref{C1} to~\ref{C4} into a Boolean ring. Recall that the notation $a+b$ is being used for $p(a,b,\bar a)$.
\begin{lemma}\label{theorem:ring}
If $(A,p,0,1)$ verifies conditions \ref{C1} to~\ref{C4} and if
\begin{equation}\label{T3-2}
p(0,a,b)=p(a,a,b)
\end{equation}
then:
\begin{equation}
\label{ring}(A,+,\wedge,0,1)\ \textrm{is a Boolean ring}.  
\end{equation}
\end{lemma}
\begin{proof}
Condition (\ref{T3-2}) implies Boolean complements and $a+a=0$:
\begin{eqnarray}
a\wedge\bar a&=&p(0,a,\bar a)\ \eq{(\ref{L3})}\ p(\bar a,\bar a,0)\ \eq{(\ref{T3-2})}\ p(0,\bar a,0)\ \eq{\ref{C2}}\ 0,\nonumber\\
\label{zero-plus}a+a&=&p(a,a,\bar a)\ \eq{(\ref{T3-2})}\ p(0,a,\bar a)=a\wedge \bar a=0.
\end{eqnarray}
Now, this result and (\ref{monoid-plus}) imply that $(a+b)+(b+a)=0$ and consequently that $+$ is commutative:
\begin{equation}
\label{commut-plus}p(a,b,\bar a)=p(b,a,\bar b).
\end{equation}
Right distributivity of $\wedge$ over $+$ is then proved as follows:
\begin{eqnarray}
(a\wedge c)+(b\wedge c)&\quad=\quad&p(a\wedge c,p(0,b,c),\overline{a\wedge c})\nonumber\\
                   &\eq{\ref{C3},\ref{C4}}&p(a\wedge c,b,p(a\wedge c,c,\overline{a\wedge c}))\nonumber\\
									 &\eq{(\ref{commut-plus})}&p(a\wedge c,b,p(c,p(0,a,c),\bar c))\nonumber\\
									 &\eq{\ref{C3},\ref{C4}}&p(a\wedge c,b,p(c,a,p(c,c,\bar c))\nonumber\\
									 &\eq{(\ref{zero-plus})}&p(p(0,a,c),b,p(c,a,0))\nonumber\\
									 &\eq{(\ref{L3})}&p(p(0,a,c),b,p(0,\bar a,c))\nonumber\\
                   &\eq{\ref{C3}}&p(0,p(a,b,\bar a),c)=(a+b)\wedge c.\label{RD}
\end{eqnarray} 
The following properties hold:
\begin{equation}
\label{LD}a\wedge(a+b)=a+(a\wedge b),\quad b\wedge(a+b)=(b\wedge a)+b.
\end{equation}
Indeed:
\begin{eqnarray*}
a\wedge(a+b)&=&p(0,a,p(a,b,\bar a)\ \eq{(\ref{zero-plus})}\ p(p(a,a,\bar a),a,p(a,b,\bar a))\\
            &\eq{\ref{C3}}&p(a,p(a,a,b),\bar a)\ \eq{(\ref{T3-2})}\ p(a,a\wedge b,\bar a)=a+(a\wedge b).
\end{eqnarray*}
The second relation in (\ref{LD}) is a consequence of the commutativity of~$+$.
We can now prove that $\wedge$ is commutative. From Lemma~\ref{lemma:complement} we already know that, under the hypothesis of Lemma~\ref{theorem:ring}, $\wedge$ is idempotent and consequently 
\begin{eqnarray*}
(a+b)\wedge (a+b)=a+b &\impli{(\ref{RD})}&\big(a\wedge(a+b)\big)+\big(b\wedge(a+b)\big)=a+b\\
                      &\impli{(\ref{LD})}&\big(a+(a\wedge b)\big)+\big((b\wedge a)+b\big)=a+b\\
                      &\impli{(\ref{monoid-plus}),(\ref{commut-plus})}&(a\wedge b)+(b\wedge a)+a+b=a+b\\
                      &\impli{(\ref{monoid-plus}),(\ref{zero-plus})}&(a\wedge b)+(b\wedge a)=0\\
											&\impli{(\ref{monoid-plus}),(\ref{zero-plus})}& a\wedge b=b\wedge a.
\end{eqnarray*}
\end{proof}

We now proceed as follows: \ref{T1} $\Rightarrow$ \ref{T2} $\Rightarrow$ (\ref{T3} $\Leftrightarrow$ \ref{T4}) $\Rightarrow$ \ref{T5} $\Rightarrow$ \ref{T1}, in order to prove Theorem \ref{thm:1}.
\begin{proof}
We begin by proving that if $(A,p,0,1)$ is a system verifying the hypothesis of Theorem \ref{thm:1} then \ref{T1} implies \ref{T2}. It is well known (see e.g. \cite{Birk-livro, NMF2012}) that, in a distributive lattice, if $x\wedge a=x'\wedge a$ and $a\vee x=a\vee x'$ for some given element $a$ in the lattice then $x=x'$. We show here that if $(A,p,0,1)$ verifies \ref{C1} to~\ref{C4} and $(A,\vee,\wedge,\bar{()},0,1)$ is a Boolean algebra then
$$\left\{\begin{array}{rcl}
p(a,b,c)\wedge c&=&((\bar{b}\wedge a)\vee(b\wedge c))\wedge c\\
c\vee p(a,b,c)&=&c\vee ((\bar{b}\wedge a)\vee(b\wedge c))
\end{array}\right.,$$
which proves $\ref{T2}$. Indeed:
$$\begin{array}{rcccl}
p(a,b,c)\wedge c&=&p(0,p(a,b,c),c)&\eq{\ref{smA7}}&p(p(0,a,c),b,p(0,c,c))\\[5pt]
                &\eq{(\ref{idem})}&p(p(0,a,c),b,c)&\eq{\ref{C4}}&p(p(0,a,c),b,p(0,1,c))\\[5pt]
								&\eq{\ref{smA7}}&p(0,p(a,b,1),c)&=&(a\vee b)\wedge c\\[5pt]
								&=&((\overline{b}\wedge a)\vee b)\wedge c
								&=&((\overline{b}\wedge a)\vee(b\wedge c))\wedge c
\end{array}$$
$$\begin{array}{rcccl}
c\vee p(a,b,c)&=&p(c,p(a,b,c),1)&\eq{\ref{smA7}}&p(p(c,a,1)),b,p(c,c,1))\\[5pt]
							&\eq{(\ref{idem})}&p(p(c,a,1),b,c)&\eq{\ref{C4}}&p(p(c,a,1),b,p(c,0,1))\\[5pt]
							&\eq{\ref{smA7}}&p(c,p(a,b,0),1)&\eq{(\ref{L3})}&p(c,p(0,\overline{b},a),1)\\[5pt]
							&=&c\vee (\overline{b}\wedge a)
							&=&c\vee((\overline{b}\wedge a)\vee (b\wedge c)).
\end{array}$$
Next, it is shown that condition $\ref{T2}$ implies condition $\ref{T3}$. Indeed, when $\ref{T2}$ is true, we have:
$$p(1,a,1)=(\bar{a}\wedge1)\vee(a\wedge 1)$$
which means, using \ref{C2} and (\ref{L7}), that $1=\bar{a}\vee a$ and, by duality, that $\bar a\e a=0$. Therefore $p(a,a,b)=(\bar a\wedge a)\vee(a\wedge b)=a\wedge b$. Conditions $\ref{T3}$ and $\ref{T4}$ are equivalent by duality (\ref{L6}). Lemma \ref{theorem:ring} proves that $\ref{T3}$ implies $\ref{T5}$. It remains to prove \ref{T5} $\Rightarrow$ \ref{T1}, that is, if $(A,+,\wedge,0,1)$ is a Boolean ring then $(A,\vee,\wedge,\overline{()},0,1)$ is a Boolean algebra with $\bar a$ defined as $a+1$ and $a\vee b$ defined as $a+b+a\wedge b$. Indeed, firstly we have:
\begin{equation}\label{a-plus-1}
a+1=p(a,1,\bar a)\ \eq{\ref{C4}}\ \bar a
\end{equation}
and consequently
\begin{equation}\label{a-plus-bar}
a+(a+1)=a+\bar a\quad\ \impli{(\ref{monoid-plus}),(\ref{zero-plus})}\quad\ a+\bar a=1.
\end{equation}
Secondly, we have
\begin{equation}\label{a-plus}
a+b+a\wedge b\ \eq{(\ref{ring})}\ a+(b\wedge(a+1))\ \eq{(\ref{a-plus-1})}\ a+(b\wedge\bar a)
\end{equation}
\begin{eqnarray*}
\Rightarrow a+b+a\wedge b&\eq{(\ref{a-plus})}& p(a,p(0,b,\bar a),\bar a)\quad\ \eq{\ref{C3},\ref{C4}}\quad\ p(a,b,p(a,\bar a,\bar a))\\
             &=& p(a,b,a+\bar a)\ \eq{(\ref{a-plus-bar})}\ p(a,b,1)=a\vee b.
\end{eqnarray*}
\end{proof}

\section{Conclusion}
It is readily shown that the category of Boolean algebras is isomorphic to the category of systems $(A,p,0,1)$ satisfying conditions \mbox{\ref{C1} to \ref{C4}} and  condition $\ref{T2}$ which can be written as 
$$p(a,b,c)=p(p(a,b,0),p(0,b,c),1),$$ 
or conditions $\ref{T3}$ or $\ref{T4}$ of Theorem \ref{thm:1}.

Other approaches to ternary Boolean algebras can be found in the literature \cite{Pad, Whiteman}. For instance, Whiteman~\cite{Whiteman} considers ternary rejection
$p(a,b,c)=\bar{a}\bar{b}+\bar{b}\bar{c}+\bar{c}\bar{a}$
which has the advantage of defining negation as $\bar{a}=p(a,a,a)$, but complete commutativity is still required. Equivalent forms of conditions \ref{B2}, \ref{B3} and \ref{B4} appear in the work of Kempe \cite{Kempe} 
and are referred as median algebras in \cite{Bandelt, Isbell, Sholander}. 

If the Boolean algebra $\ref{T1}$ in Theorem \ref{thm:1} is regarded as a Boolean ring with $ab=a\wedge b$ and $a+b=(\overline{b}\wedge a)\vee(b\wedge\overline{a})=p(a,b,\bar a)$ then 
$$(\bar{b}\wedge a)\vee(b\wedge c)=\bar{b}a+bc.$$
The expression $p(a,b,c)=\bar{b}a+bc$, with $\bar{b}=1-b$, has recently been used to model the unit interval as a mobi algebra structure~\cite{mobi} which has \ref{C1} to \ref{C4} amongst its axioms.

The ternary Mal’tsev operation~\cite{Maltsev} in a group $(ab^{-1}c)$ admits a similar characterization to the one obtained here (see Theorem 4 in \cite{Certaine}). However, as pointed out by Birkhoff and Kiss \cite{Birk}, the two operations are quite different: in a group we have $p(a,b,b)=a=p(b,b,a)$ whereas in a Boolean algebra we have $p(a,b,a)=a$. Nevertheless, as shown in~\cite{NMF2012} there are some touching points between (weakly) Mal’tsev categories and distributive lattices.


\end{document}